\UseRawInputEncoding
\documentclass[twoside,12pt]{article}
\usepackage{amssymb,amsmath,bm,euscript}
\usepackage{graphicx,color,caption}
\usepackage{ifpdf}

\usepackage[colorlinks=true]{hyperref}

\usepackage{braket}

\usepackage{algorithm, algorithmic}

\newtheorem{proposition}{Proposition}[section]
\newtheorem{theorem}[proposition]{Theorem}

\newtheorem{example}[proposition]{Example}

\newcommand{\qed}{\hphantom{.}\hfill $\Box$\medbreak}

\def\S{\mathcal{S}}

\def\I{\mathcal{I}}

\def\A{{\mathcal{A}}}
\def\B{\mathcal{B}}
\def\C{\mathcal{C}}
\def\CC{\mathbb{C}}
\def\KK{\mathbb{K}}
\def\RR{\mathbb{R}}
\def\BB{\mathbb{B}}
\def\SS{\mathbb{S}}

\def\U{\mathcal{U}}
\def\V{\mathcal{V}}

\def\S{\mathcal{S}}

\def\Y{{\mathcal{Y}}}
\def\Z{\mathcal{Z}}

\def\aa{{\bf a}}
\def\bb{{\bf b}}
\def\cc{{\bf c}}
\def\dd{{\bf d}}
\def\ii{{\bf i}}

\def\ee{{\bf e}}
\def\x{{\bf x}}
\def\y{{\bf y}}

\def\z{{\bf z}}

\def\s{{\bf s}}

\def\uu{{\bf u}}
\def\vv{{\bf v}}
\def\0{{\bf 0}}
\def\1{{\bf 1}}

\title{\bf Tubal Matrices}
\author{ \hspace{1mm}
Liqun Qi\thanks{
Department of Applied
    Mathematics, The Hong Kong Polytechnic University, Hung Hom,
    Kowloon, Hong Kong, China; ({\tt liqun.qi@polyu.edu.hk}). }
\and  and \and
Ziyan Luo\footnote{Department of Mathematics,
  Beijing Jiaotong University, Beijing 100044, China. (zyluo@bjtu.edu.cn). This author's work was supported by NSFC (Grant No.  11771038) and Beijing Natural Science Foundation (Grant No.  Z190002).}
}

\begin{document}
\date{\today}
\maketitle

\begin{abstract}
It was shown recently that the f-diagonal tensor in the T-SVD factorization must satisfy some special properties.  Such f-diagonal tensors are called s-diagonal tensors.  In this paper, we show that such a discussion can be extended to any real invertible linear transformation.   We show that two Eckart-Young like theorems hold for a third order real tensor, under any doubly real-preserving unitary transformation.  The normalized Discrete Fourier Transformation (DFT) matrix, an arbitrary orthogonal matrix, the product of the normalized DFT matrix and an arbitrary orthogonal matrix are examples of doubly real-preserving unitary transformations.  We use tubal matrices as a tool for our study.    We feel that the tubal matrix language makes this approach more natural.


\end{abstract}


\section{Introduction}

Tensor decompositions have wide applications in engineering and data science \cite{KB09}.   The most popular tensor decompositions include CP decomposition and Tucker decomposition as well as tensor train decomposition \cite{KB09, DDV00, Os11}.

The tensor-tensor product (t-product) approach, developed by Kilmer, Martin, Braman and others \cite{KMP08, Br10, KM11, KBHH13}, is somewhat different.  They defined T-product operation such that a third order tensor can be regarded as a linear operator applied on a matrix space.   They regarded vectors as  tubal scalars, matrices as tubal vectors, and third order tensors as tubal matrices, i.e., downgrade the orders of vectors, matrices and third order tensors, by one, respectively.  In \cite{KMP08, KM11}, T-SVD factorization and the best low tensor tubal rank approximation were developed.     Several important matrix terminology were extended to third order tensors, such as the transpose of a third order tensor, the inverse of a third order tensor, identity tensors, orthogonal tensors and f-diagonal tensors.
Discrete Fourier Transform (DFT) was used in this process.  After these, this approach has been widely used in applications \cite{CXZ20, HZW21, JNZH20, KKA15, LLOQ21, LYQX20, Lu18, MQW20, MQW21, QY21, QBNZ21, QBNZ21a, SHKM14, SNZ20, WCW20, XCGZ21, XCGZ21a, YHHH16, ZSKA18, ZEAHK14, ZA17, ZHW21, ZLLZ18, ZC2021}.   In 2015,
Kernfeld, Kilmer and Aeron \cite{KKA15} extended the t-product approach by defining the t-product in a so-called transform domain under an arbitrary invertible transform.   In 2020, Song, Ng and Zhang \cite{SNZ20} studied this approach by using unitary transformation matrices, while
Jiang, Ng, Zhao and Huang \cite{JNZH20} proposed to apply non-invertible transforms instead of unitary matrices.   Recently, Wang, Gu, Lee and Zhang \cite{WGLZ21} combined this approach with quantum computation.

Theoretically, this approach was vague in several places.   Somehow, in practice, vector singular values are not convenient.   Some scalar singular values of third order tensors were introduced \cite{Lu18, QY21, ZA17}.  Recently, it was shown in \cite{LLOQ21} that the f-diagonal tensor in the T-SVD factorization cannot be an arbitrary f-diagonal tensor.  The authors of \cite{LLOQ21} called such a f-diagonal tensor an s-diagonal tensor, and fully characterized such an s-diagonal tensor under the DFT.

The followings are two further theoretical issues.

First, can we characterize the f-diagonal tensor in the T-SVD factorization for the $*_L$ product developed in \cite{KKA15} with an invertible linear transformation $L$, other than the DFT?

Second, the Eckart-Young like theorem (Theorem 4.3 of \cite{KM11}) is based on the tensor tubal rank.  For an ${m \times n \times p}$ tensor $\A$, its tensor tubal rank is at most $\min \{ m , n \}$.   However, the f-diagonal tensor $\S$ in the T-SVD factorization of $\A$ has $p\min \{ m , n \}$ diagonal entries.  Is there another  Eckart-Young like theorem based upon a rank which is up to $p\min \{ m , n \}$?

We aim to answer these two questions in this paper.

We characterize the f-diagonal tensor in the T-SVD factorization of a third order real tensor for a doubly
real-preserving transformation $L$.  The DFT matrix, an arbitrary orthogonal matrix, the product of the DFT matrix and an arbitrary orthogonal matrix are examples of doubly real-preserving transformations.     Note that the Discrete Cosine Transformation (DCT) matrix is an orthogonal matrix.

We show that two Eckart-Young like theorems hold for an ${m \times n \times p}$ tensor $\A$, as described above, for any doubly real-preserving unitary transformation $L$.
Note that the real-preserving and doubly real-preserving properties may be neglected in \cite{KKA15}. This may have some potential application values.

We use tubal matrices as a tool for this discussion.  We feel that the tubal matrix language makes this approach more natural.

Throughout this paper, we assume that $m, n$ and $p$ are positive integers.  

\section{T-SVD Factorization and s-Diagonal Tubal Matrices}

\subsection{Real-Preserving Transformation and Tubal Scalar Module}

Let $L = (L_{ij}) \in \CC^{p \times p}$ be an invertible matrix.   Denote $L^{-1} \equiv H = (H_{ij})$.   For $\aa \in \CC^p$,
 define a mapping $\phi_L: \CC^p \to \CC^p$ by $\phi_L(\aa) := L\aa$, and its inverse mapping $\phi_L^{-1} : \CC^p \to \CC^p$ by $\phi_L^{-1}(\aa) := L^{-1}\aa$.   For $\aa, \bb \in \CC^p$, define
 \begin{equation} \label{e2.1}
 \aa \odot_L \bb = \phi_L^{-1}[(\phi_L(\aa)) \circ (\phi_L(\bb))],
\end{equation}
 where $\circ$ is the Hadamard product.  If $L$ is real, then $\odot_L : \RR^p \times \RR^p \to \RR^p$.   In general, even if $\aa, \bb \in \RR^p$, $\aa \odot_L \bb$ may not be real.

 For $\aa \in \CC^p$, denote its $i$th component by $\aa(i)$.

Let $F \equiv F_p$ be the $p \times p$ Discrete Fourier Transform (DFT) matrix with the form
$$F = \begin{bmatrix}
1 & 1 & \cdots & 1 & 1 \\
1 & \omega & \cdots  & \omega^{p-2} & \omega^{p-1} \\
\vdots & \vdots & \ddots& \vdots & \vdots \\
1 & \omega^{p-2} &  \cdots & \omega^{(p-2)(p-2)} & \omega^{(p-2)(p-1)}\\
1 & \omega^{p-1} &  \cdots & \omega^{(p-1)(p-2)} & \omega^{(p-1)(p-1)}
\end{bmatrix},$$
where
$$\omega = e^{-{2\pi \ii \over p}},$$
and $\ii = \sqrt{-1}$ is the imaginary unit.   Then $F^{-1} = {1 \over p}F^H$, where $F^H$ is the conjugate transpose of $F$, and for any $\aa, \cc \in \CC^p$, we have
\begin{equation} \label{e2.2}
\phi_F(\aa)(l) = \sum_{i=1}^p \omega^{(i-1)(l-1)} \aa(i),
\end{equation}
for $l=1, \cdots, p$, and
\begin{equation} \label{e2.3}
\phi_F^{-1}(\cc)(k) = {1 \over p}\sum_{l=1}^p \bar \omega^{(l-1)(k-1)} \cc(l),
\end{equation}
for $k = 1, \cdots, p$.

\begin{proposition} \label{p2.1}  For any $\aa, \bb \in \RR^p$, we have
\begin{equation} \label{e2.4}
(\aa \odot_F \bb)(k) =  \sum \left\{ \aa(i)\bb(j) : i+j-k-1= 0\ {\rm mod}(p), i, j = 1, \cdots, p \right\},
\end{equation}
for $k = 1, \cdots, p$.
Thus, $\odot_F : \RR^p \times \RR^p \to \RR^p$.
\end{proposition}
{\bf Proof}  Let $\aa, \bb \in \RR^p$, and $\cc = \phi_F(\aa) \circ \phi_F(\bb)$.   For $k = 1, \cdots, p$, we have
\begin{eqnarray*}
\phi_F^{-1}(\cc)(k) & = & {1 \over p}\sum_{l=1}^p \bar \omega^{(l-1)(k-1)} \cc(l)\\
& = & {1 \over p}\sum_{l=1}^p\sum_{i=1}^p\sum_{j=1}^p \bar \omega^{(l-1)(k-1)}\omega^{(i-1)(l-1)}\omega^{(j-1)(l-1)}\aa(i)\bb(j)\\
& = & {1 \over p}\sum_{i=1}^p\sum_{j=1}^p \aa(i)\bb(j)\left[\sum_{l=1}^p\omega^{(l-1)(i+j-k-1)}\right].
\end{eqnarray*}
Note that
$$\sum_{l=1}^p\omega^{(l-1)(i+j-k-1)} = \left\{ \begin{aligned}p, & \ \ {\rm if\ } i+j-k-1= 0, \ {\rm mod}(p) \\0, & \ \ {\rm otherwise}. \end{aligned}
\right.$$
Thus,
$$\phi_F^{-1}(\cc)(k) =  \sum \left\{ \aa(i)\bb(j) : i+j-k-1= 0\ {\rm mod}(p), i, j = 1, \cdots, p \right\}.$$
This proves (\ref{e2.4}).  The last conclusion follows.
\qed

We call a $p \times p$ invertible matrix $L$ a {\bf real-preserving transformation}  if $\odot_L : \RR^p \times \RR^p \to \RR^p$.   The DFT matrix $F$ and any $p \times p$ real invertible matrices are real-preserving transformations.    Not all complex invertible matrices are real-preserving. For example, let
$$L = \begin{bmatrix} \ii & 0 \\ 0 & \ii \end{bmatrix}, \aa = \left[
                                                                \begin{array}{c}
                                                                  1 \\
                                                                  1 \\
                                                                \end{array}
                                                              \right], \bb = \left[
                                                                \begin{array}{c}
                                                                  1 \\
                                                                  0 \\
                                                                \end{array}
                                                              \right].
$$
Then $L^{-1} = -L$, $L\aa = \ii\aa$, $L^{-1}\aa = -\ii\aa$, and
$$\aa \odot_L \bb = L^{-1}(L\aa \circ L\bb) = -\ii(-\aa \circ \bb) =  \left[
                                                                \begin{array}{c}
                                                                  \ii \\
                                                                  0 \\
                                                                \end{array}
                                                              \right],$$
which is not real.

\begin{proposition} \label{p2.2}
The product of a $p \times p$ real-preserving matrix $L$ and any $p \times p$ invertible matrix $P$ is real-preserving.
\end{proposition}
{\bf Proof}  Suppose that $\aa, \bb \in \RR^p$.   Then,
\begin{eqnarray*}
\aa \odot_{LP} \bb & = &   (LP)^{-1}\left( (LP\aa) \circ (LP\bb)\right)\\
& = & P^{-1} L^{-1}\left( (LP\aa) \circ (LP\bb)\right)\\
& = & P^{-1} \left( (P\aa) \odot_L (P\bb)\right).
\end{eqnarray*}
Since $L$ is real-preserving, $(P\aa) \odot_L (P\bb) \in \RR^p$.  Thus,
$$\aa \odot_{LP} \bb = P^{-1} \left( (P\aa) \odot_L (P\bb)\right) \in \RR^p.$$
Hence, $LP$ is real-preserving.
\qed

This shows that the set of $p \times p$ real-preserving matrices is invariant for multiplying right with any $p \times p$ invertible matrix.  In particular, the product of the DFT matrix and any invertible matrix is real-preserving. However, it is worth pointing out that $PL$ may not be real-preserving even if $P$ is an orthogonal matrix and $L$ is real-preserving, as illustrated in the following example.

\begin{example} \label{ex2.3} Let $P = \left[
                           \begin{array}{ccc}
                             \frac{\sqrt{2}}{2} & \frac{\sqrt{2}}{2} & 0\\
                             \frac{\sqrt{2}}{2} & -\frac{\sqrt{2}}{2} & 0 \\
                             0 & 0 & 1\\
                           \end{array}
                         \right]$ and $F = \left[
                                                               \begin{array}{ccc}
                                                                 1 & 1 & 1 \\
                                                                 1 & w & w^2 \\
                                                                 1 & w^2 & w^4 \\
                                                               \end{array}
                                                             \right]$. Obviously, $P$ is orthogonal and $F$ is the $3\times 3$ DFT matrix. Choose                                                              $\aa =(1, 2, 3)^{\top}$, and $\bb = (3,4,5)^{\top}$. Direct calculations lead to
\begin{equation} \aa \odot_{PF} \bb =   (PF)^{-1}\left( (PF\aa) \circ (PF\bb)\right)\approx \left[
                                                                                                      \begin{array}{c}
                                                                                                        11.4602+3.9279{\bf i} \\
                                                                                                        17.8241-8.8269{\bf i} \\
                                                                                                        22.6881+3.0619{\bf i}  \\
                                                                                                      \end{array}
                                                                                                    \right],
  \end{equation}
  which indicates that $PF$ is not real-preserving.
\end{example}

Suppose $L$ is a real-preserving transformation.   Then $\KK_p \equiv \KK_p(L) = (\RR^p, +, \odot_L)$ is a commutative ring with unity, with $+$ as the usual addition of vectors.  As discussed in \cite{Br10},  $\KK_p$ is also a module.   As in \cite{KBHH13}, we call an element $\aa \in \KK_p$ a tubal scalar, and call $\KK_p$ the $p$-dimensional tubal scalar module.  We use small bold letters $\aa, \bb, \cdots$, to denote tubal scalars.

The zero tubal scalar $\0$ has components $\0(k) = 0$ for $k = 1, \cdots, p$.
Let the tubal scalar $\1$ be defined by $\1(k) = 1$ for $k= 1, \cdots, p$.
Denote the unit tubal scalar of $\KK_p$ by $\ee_L$.

\begin{theorem}
For $k = 1, \cdots, p$, we have
\begin{equation} \label{e2.5}
\ee_L(k) = \sum_{j=1}^p H_{kj}.
\end{equation}
In particular, we have $\ee_F(1) = 1$, $\ee_F(k) = 0$, for $k = 2, \cdots, p$.
For $\aa \in \KK_p$, it is invertible if and only if $\phi_L(\aa)(j) \equiv \sum_{i=1}^p L_{ji}\aa(i) \not = 0$ for $j = 1, \cdots, p$.  In this case, we have
\begin{equation} \label{e2.6}
\aa^{-1}(k) = \sum_{j=1}^p {H_{kj} \over  \sum_{i=1}^p L_{ji}\aa(i)},
\end{equation}
for $k = 1, \cdots, p$.
\end{theorem}
{\bf Proof}  Suppose that $\ee_L$ is defined by (\ref{e2.5}).   Then $\phi_L(\ee_L) = \1$.  For any $\aa \in \KK_p$,
$$\aa \odot_L \ee_L =  \phi_L^{-1}[(\phi_L(\aa)) \circ (\phi_L(\ee_L))] = \phi_L^{-1}[(\phi_L(\aa)) \circ \1] =  \phi_L^{-1}[(\phi_L(\aa))] = \aa.$$
Similarly, $\ee_L \odot_L \aa = \aa$.  Thus, $\ee_L$ is the unit element of $\KK_p$.

Therefore, we have the formula of $\ee_L \equiv \phi_L^{-1}\left(\1\right)$.

Suppose that $\aa \in \KK_p$ is invertible.  Then
$$\aa \odot_L \aa^{-1} = \ee_L,$$
i.e.,
$$ \phi_L^{-1}\left(\phi_L(\aa) \circ \phi_L(\aa^{-1})\right) = \phi_L^{-1}\left( \1\right).$$
Thus,
$$\phi_L(\aa) \circ \phi_L(\aa^{-1}) =  \1.$$
This is equivalent to
$$\phi_L(\aa)(k)\phi_L(\aa^{-1})(k) = 1,$$
for $k=1, \cdots, p$.  The last two conclusions follow from this.
\qed

If $\aa \in \KK_p$ is not invertible, we may still derive the formula of $\aa^\dagger$, the generalized inverse of $\aa$.   We do not go to details for this.

For $\aa \in \KK_p$, define its modulus as
$$|\aa| = \sqrt{\sum_{k=1}^p \aa(k)^2}.$$

\subsection{Doubly Real-Preserving Transformation and the Transpose of a Tubal Scalar}

Suppose that $L$ is a real-preserving transformation.

For $\aa \in \KK_p$, define a mapping $\psi_L$ as
\begin{equation}
\psi_L(\aa) = \phi_L^{-1}(\overline{\phi_L(\aa)}).
\end{equation}
If $L$ is real, then $\psi_L(\aa) = \aa$.   In general, $\psi_L(\aa)$ is not real.

\begin{proposition} \label{p2.3}
The mapping $\psi_F : \KK_p \to \KK_p$.  For $\aa \in \KK_p$, we have $\psi_F(\aa)(1) = \aa(1)$, and
$\psi_F(\aa)(k)=\aa(p+2-k)$ for $k = 2, \cdots, p$.
\end{proposition}
{\bf Proof}
Let $\aa \in \KK_p$ and $\cc = \overline{\phi_F(\aa)}$.   Then for $k = 1, \cdots, p$,
\begin{eqnarray*}
\psi_F(\aa)(k) & = & \phi_F^{-1}(\cc)(k)\\
& = & {1 \over p}\sum_{l=1}^p \bar \omega^{(l-1)(k-1)} \cc(l)\\
& = & {1 \over p}\sum_{l=1}^p \bar \omega^{(l-1)(k-1)} \overline{\sum_{i=1}^p \omega^{(i-1)(l-1)} \aa(i)}\\
& = & {1 \over p}\sum_{i=1}^p\aa(i)\left[\sum_{l=1}^p \omega^{(1-l)(i+k-2)}\right].
\end{eqnarray*}
Since
$$\sum_{l=1}^p\omega^{(1-l)(i+k-2)} = \left\{ \begin{aligned}p, & \ \ {\rm if\ } i+k-2= 0, \ {\rm mod}(p) \\0, & \ \ {\rm otherwise}, \end{aligned}
\right.$$
the conclusion follows.
\qed

We say that a real-preserving transformation $L$ is a {\bf doubly real--preserving transformation} if
$\psi_L : \KK_p \to \KK_p$.   Then DFT matrix $F$ and any $p \times p$ real invertible matrices are doubly real-preserving transformations.

\begin{proposition} \label{p2.4a}
The product of a $p \times p$ doubly real-preserving matrix $L$ and any $p \times p$ invertible matrix $P$ is doubly real-preserving.
\end{proposition}
{\bf Proof}   By Proposition \ref{p2.2}, $LP$ is real-preserving.

Let $\aa \in \KK_p$.  Then,
\begin{eqnarray*}
\psi_{LP}(\aa) & = & \phi_{LP}^{-1}(\overline{\phi_{LP}(\aa)})\\
& = & P^{-1}L^{-1}(\overline{LP\aa})\\
& = & P^{-1}\left(L^{-1}(\overline{L(P\aa)})\right)\\
& = & P^{-1}\psi_L(P\aa).
\end{eqnarray*}
Since $L$ is doubly real-preserving, $\psi_L(P\aa) \in \KK_p$.  Thus,
$$\psi_{LP}(\aa) = P^{-1}\psi_L(P\aa) \in \KK_p.$$
This implies that $LP$ is doubly real-preserving.
\qed

Thus, the set of $p \times p$ doubly real-preserving matrices is invariant for multiplying right with any $p \times p$ invertible matrix.  In particular, the product of the DFT matrix and any invertible matrix is doubly real-preserving.


Suppose that $L$ is a doubly real-preserving transformation.   For any $\aa \in \KK_p$,
its transpose $\aa^\top \in \KK_p$ is defined by $\aa^\top = \psi_L(\aa)$.
If $\aa^\top = \aa$, then we say that $\aa$ is symmetric.     Thus, if $L$ is real, any $\aa \in \KK_p$ is symmetric.   The following proposition can be proved by definition.

\begin{proposition} \label{p2.4}
Suppose that $L$ is a doubly real-preserving transformation.   Let $\aa, \bb \in \KK_p$.  Then $(\aa^\top)^\top = \aa$ and $(\aa \odot_L \bb)^\top = \bb^\top \odot_L \aa^\top$.
\end{proposition}

The symmetric tubal scalars in $\KK_p$ form a subspace of $\KK_p$.   Denote it as $\SS_p$.

In general, even if $\aa \in \KK_p$, $\phi_L(\aa)$ is not real.   However, we have the following theorem.

\begin{theorem}
Suppose that $L$ is a doubly real-preserving transformation.
The mapping $\phi_L$ defined above maps $\SS_p$ to $\KK_p$.   Furthermore, we have
$\phi_L : \SS_p \to \SS_p$ if $L$ is real or $L =F$.
\end{theorem}
{\bf Proof}  Let $\aa \in \SS_p$.  Then
$$\phi_L(\aa) = \phi_L(\aa^\top)= \phi_L\left(\phi_L^{-1}(\overline{\phi_L(\aa)})\right) = \overline{\phi_L(\aa)}.$$
Thus, $\phi_L(\aa)$ is real.  If $L$ is real, then $\KK_p(L) \equiv \SS_p(L)$.
We also have
\begin{equation*}
\phi_F(\aa)(1)= \sum_{l=1}^p \omega^{(l-1)(1-1)}\aa(l)= \sum_{l=1}^p \aa(l),
\end{equation*}
and for $k = 2, \cdots, p$,
\begin{eqnarray*}
\phi_F(\aa)(k)& = & \aa(1) + {1 \over 2}\sum_{l=2}^p \left[\omega^{(l-1)(k-1)} +  \bar \omega^{(l-1)(k-1)}\right]\aa(l)\\
& = & \phi_F(\aa)(p-k+2).
\end{eqnarray*}
Thus, $\phi_F(\aa) \in \SS_p$.
\qed

We say that $\aa \in \SS_p$ is positive semi-definite if $\phi_L(\aa)(k) \ge 0$ for $k = 1, \cdots, p$.

\begin{theorem}
Suppose that $L \in \CC^{p \times p}$ is a doubly real-preserving transformation, and
$\aa \in \KK_p$.   Then $\aa^\top \odot_L \aa$ is symmetric and positive semi-definite.
\end{theorem}
{\bf Proof}
 By Proposition \ref{p2.4},
 $$(\aa^\top \odot_L \aa)^\top = \aa^\top \odot_L (\aa^\top)^\top = \aa^\top \odot_L \aa.$$
 Hence, $\aa^\top \odot_L \aa$ is symmetric.  Also,
\begin{eqnarray*}
 \phi_L(\aa^\top \odot_L \aa)
& = & \phi_L\left( \phi_L^{-1}[\phi_L(\aa^\top) \circ \phi_L(\aa)]\right)\\
& = & \left[\phi_L(\aa^\top) \circ \phi_L(\aa)\right]\\
& = & \left[\phi_L(\psi_L(\aa)) \circ \phi_L(\aa)\right]\\
& = & \left[\phi_L(\phi_L^{-1}(\overline{\phi_L(\aa)}) \circ \phi_L(\aa)\right]\\
& = & \left[\overline{\phi_L(\aa)} \circ \phi_L(\aa)\right]\\
& \ge & \0.
\end{eqnarray*}
This shows that $\aa^\top \odot_L \aa$ is positive semi-definite.
\qed

\subsection{Tubal Vectors}

Suppose that $L$ is a doubly real-preserving transformation.

As in \cite{KBHH13}, an element $X \in \KK_p^n \equiv \KK_p^n(L)$ is an $n$-dimensional tubal vector with $p$-dimensional tubal scalar as its components:
$$X = \begin{bmatrix}
\x_1\\
\x_2\\
\vdots\\
\x_n\\
\end{bmatrix}.$$
We call $X$ an $n$-dimensional column tubal vector, and
$$X^\top = (\x_1^\top, \x_2^\top, \cdots, \x_n^\top)$$
an $n$-dimensional row tubal vector.  Then $\KK_p^n$ is a linear vector space for addition and multiplication with real numbers.  We use capital letters $X, Y, Z, \cdots$, to denote tubal vectors.

For $X \in \KK_p^n$, if
$$X^\top *_L X \equiv \sum_{k=1}^n \x_k^\top \odot_L \x_k = \ee_L,$$
then we say that $X$ is normalized.

For $X, Y \in \KK_p^n$, if
$$X^\top * Y \equiv \sum_{k=1}^n \x_k^\top \odot \y_k = \0,$$
then we say that $X$ and $Y$ are orthogonal.


For $X \in \KK_p^n$, define
$${\rm unfold}(X) = \begin{bmatrix}
\x^{(1)}\\
\x^{(2)}\\
\vdots\\
\x^{(p)}\\
\end{bmatrix},$$
where
$$\x^{(k)}  = \begin{bmatrix}
\x_1(k)\\
\x_2(k)\\
\vdots\\
\x_n(k)\\
\end{bmatrix} \in \RR^n.$$
Thus, ${\rm unfold}(X) \in \RR^{np}$.  Note here $\x^{(k)}$ is a vector in $\RR^n$, not a tubal scalar.  Then we define $X = {\rm fold}({\rm unfold}(X))$.

\subsection{Tubal Matrices}

Suppose that $L$ is a doubly real-preserving transformation.

As in \cite{KBHH13}, an element $\A \in \KK_p^{m \times n} \equiv \KK_p^{m \times n}(L)$ is an $(m \times n)$-dimensional tubal matrix with $p$-dimensional tubal scalar as its components:
$$\A = (\aa_{ij}) = \begin{bmatrix}
\aa_{11} & \aa_{12} & \cdots & \aa_{1n}\\
\aa_{21} & \aa_{22} & \cdots & \aa_{2n}\\
\vdots & \vdots & & \vdots \\
\aa_{m1} & \aa_{m2} & \cdots & \aa_{mn}\\
\end{bmatrix}.$$
We say that
$$A_i = (\aa_{i1}, \aa_{i2}, \cdots, \aa_{in})$$
is the $i$th row tubal vector of $\A$.   Similarly, we may define the $j$th colomn tubal vector of $\A$.

Here, we abuse the notation somewhat as follows.   We also refer $\A = (a_{ijk})$ as a third order tensor in $\RR^{m \times n \times p}$ with $a_{ijk} \equiv \aa_{ij}(k)$.  We will see that the operations and subclasses of tubal matrices are corresponding to the $*_L$-product operations and subclasses of third order tensors.

We use calligraphic letters $\A, \B, \U, \V, \S, \cdots$ to denote tubal matrices.

The addition of two tubal matrices, and the multiplication of a tubal matrix and a real number are defined directly.   Let $\A = (\aa_{ij}) \in \KK_p^{m \times n}$ and $X = (\x_j) \in \KK_p^n$.  Then their multiplication is defined as
$$\A *_L X = \left(\sum_{j=1}^n \aa_{ij} \odot_L \x_j\right) \in \KK_p^m.$$
Let $\A = (\aa_{il}) \in \KK_p^{m \times s}$ and $\B = (\bb_{lj}) \in \KK_p^{s \times n}$.   Then $\A *_L \B = \C \equiv (\cc_{ij}) \in \KK_p^{m \times n}$ be defined by
$$\cc_{ij} = \sum_{l=1}^s \aa_{il} \odot_L \bb_{lj},$$
for $i = 1, \cdots, m$ and $j = 1, \cdots, n$.    

Let $\A = (\aa_{ij}) \in \KK_p^{m \times n}$.    We call $\aa_{ii}$, for $i = 1, \cdots, \min \{ m, n \}$, diagonal entries of $\A$, and $\aa_{ij}$, for $i \not = j$, off-diagonal entries of $\A$.
We call $\A = (\aa_{ij}) \in \KK_p^{m \times n}$ a diagonal tubal matrix, if all of its off-diagonal entries are $\0$.   We may also denote a diagonal tubal matrix $\A = (\aa_{ij}) \in \KK_p^{m \times n}$ by $\A = {\rm diag}( \aa_{11}, \cdots, \aa_{\min\{m,n\}, \min\{m,n\}})$.

We call $\A \in \KK_p^{n \times n}$ an identity tubal matrix, if it is diagonal and all of its diagonal entries are $\ee_L$. We denote it as $\I_n \equiv \I_n(L)$.
Let $\A, \B \in \KK_p^{n \times n}$.  If $\A *_L \B = \B *_L \A = \I_n$, then we say that $\B$ is the inverse of $\A$, and denote $\A^{-1} = \B$.

Let $\A = (\aa_{ij}) \in \KK_p^{m \times n}$.  Then $\B = (\bb_{ij}) \in \KK_p^{n \times m}$ is defined as the transpose of $\A$ and denoted as $\A^\top = \B$, if $\bb_{ij} = \aa_{ji}^\top$ for $i = 1, \cdots, n$ and $j= 1, \cdots, m$.

If $\A = (\aa_{ij}) \in \KK_p^{n \times n}$ and $\A^{-1} = \A^\top$, then $\A$ is called an orthogonal tubal matrix.   If $\A = (\aa_{ij}) \in \KK_p^{n \times n}$ and $\A = \A^\top$, then $\A$ is called a symmetric tubal matrix.

If $L=F$, then the definitions of diagonal tubal matrices, identity tubal matrices, the inverse and the transpose of a tubal matrix, and orthogonal tubal matrices, are corresponding to the definitions of f-diagonal tensors, identity tensors, the inverse and the transpose of a third order tensor, and orthogonal tensors in \cite{Br10, KM11, KBHH13}.

\begin{proposition} \label{p2.7}
Suppose that $L$ is a doubly real-preserving transformation.
Let $\A \in \KK_p^{n \times n}$.   Then $\A$ is an orthogonal tubal matrix if and only if its row (column) tubal vectors are normalized and orthogonal to each other.
\end{proposition}
{\bf Proof} This holds directly from the multiplication rules of tubal matrices.
\qed

\begin{proposition} \label{p2.8}
Suppose that $L$ is a doubly real-preserving transformation.
Then $\A = (\aa_{ij}) \in \KK_p^{n \times n}$ is an orthogonal tubal matrix, if and only if for $k = 1, \cdots, p$,
$(\phi_L(\aa_{ij}(k))$ are $n \times n$ unitary matrices.
\end{proposition}
{\bf Proof} Suppose that $\A$ is an orthogonal tubal matrix.   Let $\B \equiv \A^\top = (\bb_{jl})$.
Then for $i, l = 1, \cdots, n$,
$$\sum_{j=1}^n \aa_{ij} \odot_L \bb_{jl} = \sum_{j=1}^n \phi_L^{-1}(\phi_L(\aa_{ij}) \circ \phi_L(\bb_{jl})) = \delta_{il}\ee_L,$$
where $\delta_{il}$ is the Kronecker symbol.    Since $\phi_L$ is linear, this is equivalent to
$$\sum_{j=1}^n  \phi_L(\aa_{ij})(k) \phi_L(\bb_{jl})(k) = \delta_{il},$$
for $k = 1, \cdots, p$, $i, l = 1, \cdots, n$, i.e.,  $(\phi_L(\aa_{ij}(k))$ are $n \times n$ unitary matrices, for $k = 1, \cdots, p$.

On the other hand, suppose that  $(\phi_L(\aa_{ij}(k))$ are $n \times n$ unitary matrices, for $k = 1, \cdots, p$.   Suppose that $(\phi_L(\aa_{ij}(k))^{H} = (\vv_{jl}(k))$, where $(\cdot)^{H}$ is the conjugate transpose, for $k = 1, \cdots, p$.   Then
$$\sum_{j=1}^n \phi_L(\aa_{ij})(k)\vv_{jl}(k) = \delta_{il},$$
for $k = 1, \cdots, p$, $i, l = 1, \cdots, n$.   These imply that for $i, l = 1, \cdots, n$,
$$\sum_{j=1}^n \aa_{ij} \odot_L \phi_L^{-1}(\vv_{jl}) = \sum_{j=1}^n \phi_L^{-1}(\phi_L(\aa_{ij}) \circ \vv_{jl}) = \delta_{il}\ee_L,$$
i.e., $\A$ is an orthogonal tubal matrix.
\qed

We now define s-diagonal tubal matrices.     Let $\A = (\aa_{ij}) \in \KK_p^{m \times n}$ be a diagonal tubal matrix.    We say that $\A$ is an s-diagonal tubal matrix if it satisfies the following properties:
\begin{itemize}
\item{(1)} All of the diagonal entries of $\A$ are symmetric;

\item{(2)} All of the diagonal entries of $\A$ are positive semi-definite;

\item{(3)} For $k=1, \cdots, p$, we have
$$\phi_L(\aa_{11})(k) \ge \phi_L(\aa_{22})(k) \ge \cdots \ge \phi_L(\aa_{\min\{m, n\}, \min\{m, n\}})(k).$$

\end{itemize}

The identity tubal matrices are s-diagonal tubal matrices.  Furthermore, in the case that $L=F$, these three properties, if translated to the language of third order tensors, are exactly the three sufficient and necessary conditions of s-diagonal tensors, identified in \cite[Theorem 6.1]{LLOQ21}.  Thus, in the case that $L=F$, an s-diagonal tubal matrix is corresponding to an s-diagonal tensor defined in \cite{LLOQ21}, and vice versa.

\subsection{T-SVD Factorization and s-Diagonal Tubal Matrices}

Suppose that $L$ is a doubly real-preserving transformation.  Let $\A = (\aa_{ij}) \in \KK_p^{m \times n}$.   Let
\begin{equation} \label{e2.9}
\Phi_L(\A) = (\phi_L(\aa_{ij})).
\end{equation}
Denote
\begin{equation} \label{e2.10}
\Phi_L(\A)(k) = (\phi_L(\aa_{ij})(k)),
\end{equation}
for $k = 1, \cdots, p$. Then the SVD of $\Phi_L(\A)(k) \in \CC^{m \times n}$ has the form
\begin{equation} \label{e2.11}
\Phi_L(\A)(k) = \Xi(k)D(k)\Theta(k)^H,
\end{equation}
where $\Xi(k) \in \CC^{m \times m}$, $\Theta(k) \in \CC^{n \times n}$ are unitary matrices,
$$D(k) = {\rm diag}(\dd_{11}(k), \cdots, \dd_{\min\{m,n\},\min\{m,n\}}(k))$$
are $m \times n$ diagonal nonnegative matrices, and
\begin{equation} \label{e2.12}
\dd_{11}(k) \ge \cdots \ge \dd_{\min\{m,n\},\min\{m,n\}}(k) \ge 0, ~~\forall k = 1, \cdots, p.
\end{equation}
Let $\S = {\rm diag}(\s_{11}, \cdots, \s_{\min\{m,n\},\min\{m,n\}})$, where $\s_{ll} = \phi_L^{-1}(\dd_{ll})$ for $l = 1, \cdots, \min \{m, n\}$.   Then we see that $\S$ is an s-diagonal tubal matrix.   Let $\Xi(k) = (\Xi_{il}(k))$, for $k = 1, \cdots, p$, $\uu_{il} = \phi_L^{-1}(\Xi_{il})$ and $\U = (\uu_{il}) \in \KK_p^{m \times m}$;  $\Theta(k) = (\Theta_{lj}(k))$, for $k = 1, \cdots, p$, $\vv_{lj} = \phi_L^{-1}(\Theta_{lj})$ and $\V = (\vv_{lj}) \in \KK_p^{n \times n}$  By Proposition \ref{p2.8}, $\U$ and $\V$ are orthogonal tubal matrices.
Let
$$\B = (\bb_{ij})= \U *_L \S *_L \V^\top \in \KK_p^{m \times n}.$$

Then
\begin{eqnarray*}
\bb_{ij} & = & \sum_{l=1}^{\min \{m, n\}} (\uu_{il} \odot_L \s_{ll}) \odot_L \vv_{jl}^\top\\
& = & \sum_{l=1}^{\min \{m, n\}} \left(\phi_L^{-1}(\phi_L(\uu_{il}) \circ \phi_L(\s_{ll}))\right) \odot_L \vv_{jl}^\top\\
& = & \sum_{l=1}^{\min \{m, n\}}  \left(\phi_L^{-1}\left(  \Xi_{il} \circ \dd_{ll}\right)\right) \odot_L \psi_L(\vv_{jl})\\
& = & \sum_{l=1}^{\min \{m, n\}} \phi_L^{-1}\left((\Xi_{il} \circ \dd_{ll}) \circ \overline{\phi_L(\vv_{jl})}\right)\\
& = & \sum_{l=1}^{\min \{m, n\}} \phi_L^{-1}\left(\Xi_{il} \circ \dd_{ll} \circ \overline{\Theta_{jl}}\right)\\
& = &  \phi_L^{-1}\left(\sum_{l=1}^{\min \{m, n\}}\Xi_{il} \circ \dd_{ll} \circ \overline{\Theta_{jl}}\right).\\
\end{eqnarray*}
Hence, for $k = 1, \cdots, p$,
$$\Phi_L(\U *_L \S *_L \V^\top)(k) = \left(\phi_L(\bb_{ij})(k)\right) = \left(\sum_{l=1}^{\min \{m, n\}} \Xi_{il}(k)\dd_{ll}(k)\overline{\Theta_{jl}(k)}\right) = \Xi(k)D(k)\Theta(k)^H =\Phi_L(\A)(k).$$
Thus,
\begin{equation} \label{e2.13}
\A = \U *_L \S *_L \V^\top \in \KK_p^{m \times n}.
\end{equation}

This is the T-SVD factorization of $\A$, under the transformation $L$.

Denote $\S = G(\A)$.   Let $\A, \B \in \KK_p^{m \times n}$.   We say that $\A$ and $\B$ are orthogonally equivalent if there are orthogonal tubal matrices $\Y \in \KK_p^{m \times m}$ and $\Z \in \KK_p^{n \times n}$ such that $\A = \Y *_L \B *_L \Z^\top$.  We first extend Theorem 3.2 of \cite{QLLO21} to doubly real-preserving transformations.

\begin{theorem} \label{t2.9}
Suppose that $L$ is a doubly real-preserving transformation.   If $\A, \B \in \KK_p^{m \times n}$ are orthogonally equivalent, then $G(\A) = G(\B)$.
\end{theorem}
{\bf Proof}   Suppose that $\A = \Y *_L \B *_L \Z^\top$, where  $\Y = (\y_{ij}) \in \KK_p^{m \times m}$ and $\Z = (\z_{ij}) \in \KK_p^{n \times n}$ are orthogonal tubal matrices.    Then we have
$$(\phi_L(\aa_{ij})(k)) = (\phi_L(\y_{ij})(k))(\phi_L(\bb_{ij})(k))(\phi_L(\z_{ij})(k))^{H},$$
where by Proposition \ref{p2.8}, $(\phi_L(\y_{ij})(k))$ and $(\phi_L(\z_{ij})(k))$ are unitary matrices, for $k = 1, \cdots, p$.   This implies that for $k = 1, \cdots, p$, $(\phi_L(\aa_{ij})(k))$ and $(\phi_L(\bb_{ij})(k))$ have the same set of singular values, i.e., for (\ref{e2.11}), the $D$ part is invariant if $\A$ is replaced by $\B$.   Thus,  $G(\A) = G(\B)$.
\qed

\begin{theorem} \label{t2.10}
Suppose that $L$ is a doubly real-preserving transformation, and $\A \in \KK_p^{m \times n}$ has T-SVD factorization (\ref{e2.13}).  Then $\S = G(\S)$.   If $L=F$ or $L$ is real, then $\S$ is an s-diagonal tubal matrix.
\end{theorem}
{\bf Proof}  By (\ref{e2.13}), $\S$ and $\A$ are orthogonally equivalent.  By Theorem \ref{t2.9}, $\S \equiv G(\A) = G(\S)$.  If $L=F$, then $\S$ is an s-diagonal tubal matrix by Theorem 6.1 of \cite{LLOQ21}.
If $L$ is real, then all tubal scalars in $\KK_p$ are symmetric.   By (\ref{e2.12}) and the definition of s-diagonal tubal matrix, $\S$ is an s-diagonal tubal matrix.
\qed

\section{Eckart-Young Like Theorems}

\subsection{Doubly Real-Preserving Unitary Transformation}

Suppose that $L \in \CC^{p \times p}$ is doubly real-preserving and unitary.   We say that $L$ is a doubly real-preserving unitary transformation.  The normalized DFT matrix ${1\over \sqrt{p}} F$, the discrete cosine transformation (DCT) matrix, and any orthogonal matrices are doubly real-preserving.  We have the following proposition.

\begin{proposition} \label{p3.1}
The product of a $p \times p$ doubly real-preserving unitary matrix $L$ and any $p \times p$ orthogonal matrix $P$ is a doubly real-preserving unitary matrix.
\end{proposition}
{\bf Proof}  By Proposition \ref{p2.4a}, $LP$ is doubly real-preserving.  Since both $L$ and $P$ are unitary, $LP$ is also unitary.
\qed

Therefore, the set of $p \times p$ doubly real-preserving unitary matrices is invariant for multiplying right with any $p \times p$ orthogonal matrix.  In particular, the product of the DFT matrix and any orthogonal matrix is doubly real-preserving and unitary.

For $\A = (\aa_{ij}) \in \KK_p^{m \times n}$, the Frobenius norm of $\A$ is defined to be $\|\A\|_F = \sqrt{\sum_{i=1}^m\sum_{j=1}^n |\aa_{ij}|^2}$.    We have the following propositions.

\begin{proposition}
Suppose that $L$ is a doubly real-preserving unitary transformation.

(a) If $\aa, \bb \in \KK_p$, then $|\aa \odot_L \bb| \le |\aa||\bb|$.

(b) If $\A \in \KK_p^{m \times r}$ and $\B \in \KK_p^{r \times n}$, then $\|\A *_L \B\|_F \le \|\A\|_F\|\B\|_F$.

(c) If $\A \in \KK_p^{m \times m}$ is an orthogonal tubal matrix, and $\B \in \KK_p^{m \times n}$, then
$\|\A *_L \B\|_F = \|\B\|_F$.
\end{proposition}
{\bf Proof}
(a) Since $L$ is unitary,
$$|\aa \odot_L \bb| = |\phi_L^{-1}(\phi_L(\aa) \circ \phi_L(\bb))| = |\phi_L(\aa) \circ \phi_L(\bb)|
\le |\phi_L(\aa)||\phi_L(\bb)| = |\aa||\bb|.$$
(b) and (c) can be proved like in matrix theory.
\qed

\begin{proposition}
Suppose that $L$ is a doubly real-preserving unitary transformation, and $\S = (\s_{ij}) \in \KK_p^{m \times n}$ is an s-diagonal tubal matrix.   Then $\S$ has the following {\bf Tubal Scalar Value Decay Property}:
\begin{equation} \label{e3.14}
|\s_{11}| \ge |\s_{22}| \ge \cdots \ge |\s_{\min\{m,n\},\min\{m,n\}}|.
\end{equation}
\end{proposition}
{\bf Proof}
Since $\S$ is an s-diagonal tubal matrix, there is $\A \in \KK_p^{m \times n}$ such that $\S = G(\A)$. Then we have $D(k)$ for $k = 1, \cdots, p$, such that (\ref{e2.12}) holds.   This implies that
$$|\dd_{11}| \ge |\dd_{22}| \ge \cdots \ge |\dd_{\min\{m,n\},\min\{m,n\}}|,$$
where $\s_{ii} = \phi_L^{-1}(\dd_{ii})$ for $i = 1, \cdots, \min \{m, n\}$.
Since $L$ is unitary, we have $|\dd_{ii}| = |\s_{ii}|$ for $i = 1, \cdots, \min \{m, n\}$.   Then
(\ref{e3.14}) follows.
\qed

For $\A \in \KK_p^{m \times n}$, suppose $\S = G(\A)$.   Then $\S$ is an s-diagonal tubal matrix.   Let the diagonal entries of $\S$ be $\s_{11}, \cdots, \s_{\min\{m,n\}, \min\{m,n\}}$. Then the $i$th largest {\bf T-singular value} of $\A$ is defined as $\sigma_i = |\s_{ii}|$, for $i = 1, \cdots, \min \{m,n\}$.  By the tubal scalar value decay property (\ref{e3.14}), we have
\begin{equation}
\sigma_1 \ge \sigma_2 \ge \cdots \ge \sigma_{\min\{m,n\}}.
\end{equation}
The tensor tubal rank of a third order tensor  is a critical concept for the best low tensor rank theory built in \cite{KM11}.   The tensor rank in Theorem 4.3 of \cite{KM11} is called the tensor tubal rank in the later papers.  Correspondingly, we may have the tubal rank of a tubal matrix.
The {\bf tubal rank} of $\A$ is equal to the number of its nonzero T-singular values.  The $i$th tail energy of $\A$ is defined as $\tau_i = \sqrt{\sum_{j=i+1}^{\min\{m,n\}} \sigma_j^2}$.  Then the low tensor tubal rank approximation theorem, Theorem 4.3 of \cite{KM11} - the First Eckart-Young like theorem for third order tensors, can be stated in a way very similar to its matrix counterpart - Theorem 2.4.8 of \cite{GV13}.  We do this in the next subsection.

\subsection{Tubal Rank and The First Eckart-Young Like Theorem}

Suppose that $L$ is a doubly real-preserving unitary transformation, and $\A \in \KK_p^{m \times n}$. Denote the tubal rank of $\A$ by ${\rm Rank}_t(\A)$.  Assume that $\Phi_L(\A)$ and $\Phi_L(\A)(k)$ for $k = 1, \cdots, p$ are defined by (\ref{e2.9}) and (\ref{e2.10}).  By the
$${\rm Rank}_t(\A) = \max \{ {\rm Rank}(\Phi_L(\A)(k)) : k = 1, \cdots, p \}.$$

\begin{theorem}
Suppose that $L$ is a doubly real-preserving unitary transformation.   Then we have the following three conclusions.

(a) For any $\A \in \KK_p^{m \times n}$, if ${\rm Rank}_t(\A) = r$, then there are $\B \in \KK_p^{m \times r}$ and $\C \in \KK_p^{r \times n}$ such that $\A = \B *_L \C$ and ${\rm Rank}_t(\B) = {\rm Rank}_t(\C) = r$.

(b) For any $\B \in \KK_p^{m \times r}$ and $\C \in \KK_p^{r \times n}$, let $\A = \B *_L \C$.  Then
$${\rm Rank}_t(\A) \le \min \{ {\rm Rank}_t(\B), {\rm Rank}_t(\C) \}.$$

(c) For any $\A \in \KK_p^{m \times n}$,
$${\rm Rank}_t(\A) = \min \{ r : \A = \B *_L \C, \B \in \KK_p^{m \times r}, \C \in \KK_p^{r \times n} \}.$$
\end{theorem}
{\bf Proof} (a) Denote ${\rm Rank}(\Phi_L(\A)(k)) = r_k$.  Then $r = \max \{ r_k : k = 1, \cdots, p \}$.
Since  ${\rm Rank}(\Phi_L(\A)(k)) = r_k$, the matrix $(\Phi_L(\A)(k)) = \bar B_k \bar C_k$, where $\bar B_k \in \CC^{m \times r_k}$, $\bar C_k \in \CC^{r_k \times n}$, and ${\rm Rank}({\bar B_k}) = {\rm Rank}({\bar C_k}) = r_k$.  We further have $(\Phi_L(\A)(k)) = \hat B_k \hat C_k$, where $\hat B_k = (\bar B_k, O) \in \CC^{m \times r}$, $\hat C_k = ({\bar C_k \atop O}) \in \CC^{r \times n}$, and
$$\max \{ {\rm Rank}(\hat B_k) =  {\rm Rank}(\hat C_k) : k = 1, \cdots, p \} = r.$$
Define $\hat \B = (\hat \bb_{il})$ by $\hat \bb_{il}(k) = \hat B_k(i,l)$, and $\hat \C = (\hat \cc_{lj})$ by $\hat \cc_{lj}(k) = \hat C_k(l,j)$,
$\B = (\bb_{il}) \in \KK_p^{m \times r}$ by $\bb_{il} = \phi_L^{-1}(\hat \bb_{il})$, and $\C = (\cc_{lj}) \in \KK_p^{r \times n}$ by $\cc_{lj} = \phi_L^{-1}(\hat \cc_{lj})$.  Then  $\A = \B *_L \C$,
$${\rm Rank}_t(\B) = \max \{ {\rm Rank}(\Phi_L(\B)(k)) : k = 1, \cdots, p \} = \max \{ {\rm Rank}(\hat B_k) : k = 1, \cdots, p \} = r,$$
$${\rm Rank}_t(\C) = \max \{ {\rm Rank}(\Phi_L(\C)(k)) : k = 1, \cdots, p \} = \max \{ {\rm Rank}(\hat C_k) : k = 1, \cdots, p \} = r.$$

(b) We have
$${\rm Rank}_t(\A) = \max \{ {\rm Rank}(\Phi_L(\A)(k)) : k = 1, \cdots, p \},$$
$${\rm Rank}_t(\B) = \max \{ {\rm Rank}(\Phi_L(\B)(k)) : k = 1, \cdots, p \},$$
$${\rm Rank}_t(\C) = \max \{ {\rm Rank}(\Phi_L(\C)(k)) : k = 1, \cdots, p \},$$
$\Phi_L(\A)(k) = \Phi_L(\B)(k)\Phi_L(\C)(k)$,
$${\rm Rank}(\Phi_L(\A)(k)) \le \min \{ {\rm Rank}(\Phi_L(\B)(k)), {\rm Rank}(\Phi_L(\C)(k)) \},$$
for $k = 1, \cdots, p$.   The
conclusion follows.

(c) follows from (a) and (b).
\qed

From this, we have the first Eckart-Young like theorem for $*_L$ product.

\begin{theorem} \label{t3.4}
Suppose that $L$ is a doubly real-preserving unitary transformation, and
$\A \in \KK_p^{m, n}$ has a T-SVD factorization (\ref{e2.13}), where
\begin{equation}\label{S_D}\S = {\rm diag}( \s_{11}, \cdots, \s_{\min\{m,n\}, \min\{m,n\}}) = G(\A)
\end{equation} is an s-diagonal tubal matrix.  For $1 \le i < \min\{m,n\}$, let $\S_i = {\rm diag}( \s_{11}, \cdots, \s_{ii}, 0, \cdots, 0)$ and $\A_i = \U * \S_i * \V^\top$.   Then $\A_i$ is the best tubal rank $i$ approximation of $\A$, with
$$\|\A - \A_i\|_F^2 = \tau_i^2.$$
\end{theorem}

The proof of this theorem is similar to the proof of Theorem 4.3 of \cite{KM11}.  We do not go to details.
Note that in \cite{KM11}, the DFT matrix is not a unitary matrix, but a unitary matrix multiplied with a positive constant.   But it is not difficult to extend the discussion in this section to doubly real-preserving transformation multiplied with a positive constant.

\subsection{B-Rank and The Second Eckart-Young Like Theorem}

For $\A \in \KK_p^{m \times n}$, suppose $\S = G(\A)$.   Then $\S$ is an s-diagonal tubal matrix.   Let the diagonal entries of $\S$ be $\s_{11}, \cdots, \s_{\min\{m,n\}, \min\{m,n\}}$, and $\dd_{ii} = \phi_L(\s_{ii})$, for $i = 1, \cdots, \min \{m,n\}$.   Then the $j$th largest {\bf B-singular value} $\mu_j$ of $\A$ is defined as $j$th largest number of $\{ \dd_{11}(1), \cdots, \dd_{11}(p), \dd_{22}(1), \cdots,
\dd_{\min\{m,n\}, \min\{m,n\}}(p) \}$.   The {\bf B-rank} of $\A$ is equal to the number of its nonzero B-singular values, and denoted as ${\rm Rank}_b(\A)$.   The $i$th B-tail energy of $\A$ is defined as $\nu_i = \sqrt{\sum_{j=i+1}^{p\min\{m,n\}} \mu_j^2}$.   Here, the letter B means block diagonalization.
It involves a block diagonal matrix below.

Denote $J_1 = \{ 1, \cdots, \min \{ m, n \} \}, J_2 = \{ 1, \cdots, p \}$ and $J = \{ 1, \cdots, p\min \{ m, n \} \}$. Now, let $\eta : J_1 \times J_2 \to J$ be a function such that $\dd_{ii}(k) = \mu_{\eta(i, k)}$.  Such a function $\eta$ may not be unique, but we fix one function $\eta$.

 Let $\A \in \KK_p^{m \times n}$.   Define ${\rm Bldg}(\Phi_L(\A))$ as a
block diagonal matrix
$${\rm Bldg}(\Phi_L(\A)) := {\rm diag}(\Phi_L(\A)(1), \cdots, \Phi_L(\A)(p)).$$

Then, the B-singular values of $\A$ are the singular values of the matrix ${\rm Bldg}(\Phi_L(\A))$, the B-rank of $\A$ is the rank of ${\rm Bldg}(\Phi_L(\A))$.   Since $L$ is unitary, we have
$$\|\A\|_F = \|\Phi_L(\A)\|_F = \|{\rm Bldg}(\Phi_L(\A))\|_F.$$

Also, if $\A$ has T-SVD factorization (\ref{e2.13}), we have $\|\A\|_F = \|\S\|_F$.

Based upon these and the Eckart-Young theorem of the matrix ${\rm Bldg}(\Phi_L(\A))$ \cite{GV13},
we have the second Eckart-Young like theorem for $*_L$ products as follows.

\begin{theorem} \label{t3.5}
Suppose that $L$ is a doubly real-preserving unitary transformation, and
$\A \in \KK_p^{m, n}$ has a T-SVD factorization (\ref{e2.13}), where  $\S$ as defined in \eqref{S_D} 
is an s-diagonal tubal matrix.  For $1 \le j < p\min\{m,n\}$, let $\dd_{ii}^{(j)} \in \KK_p$ be defined as $\dd_{ii}^{(j)}(k) = \dd_{ii}(k)$ if $\eta(i, k) \le j$; otherwise we have $\dd_{ii}^{(j)}(k) = 0$.   Let $\s_{ii}^{(j)} = \phi_L^{-1}(\dd_{ii}^{(j)})$ for $i = 1, \cdots, \min \{m, n \}$,
$\S_j = {\rm diag}\left( \s_{11}^{(j)}, \cdots, \s_{\min\{m,n\}, \min\{m,n\}}^{(j)}\right)$ and $\A_j = \U * \S_j * \V^\top$.   Then $\A_j$ is the best B-rank $j$ approximation of $\A$, with
$$\|\A - \A_j\|_F^2 = \nu_j^2.$$
\end{theorem}
{\bf Proof} Denote the set of all tubal matrices in $\KK_p^{m \times n}$ with B-rank no more than $j$ by $\BB_j$.
Let
$$\Xi := {\rm diag}(\Xi(1), \cdots, \Xi(p)) \in \CC^{mp \times mp},$$
$$D := {\rm diag}(D(1), \cdots, D(p)) \in \RR^{mp \times np}$$
and
$$\Theta := {\rm diag}(\Theta(1), \cdots, \Theta(p)) \in \CC^{np \times np},$$
where $\Xi(k), D(k)$ and $\Theta(k)$ for $k = 1, \cdots, p$, are defined by (\ref{e2.11}).
It is easy to verify that $\Xi$ and $\Theta$ are unitary matrices, and
$${\rm Bldg}(\Phi_L(\A)) = \Xi D \Theta^H.$$
Direct manipulation yields
\begin{eqnarray*}
&& \min_{\B \in \BB_j} \|\A-\B\|_F^2\\
 & = & \min_{\B \in \BB_j} \|{\rm Bldg}(\Phi_L(\A))- {\rm Bldg}(\Phi_L(\B))\|_F^2\\
 & = & \min_{\B \in \BB_j} \|\Xi D \Theta^H- {\rm Bldg}(\Phi_L(\B))\|_F^2\\
  & = & \min\limits_{\|(\hat \dd_{11}(1), \cdots, \hat \dd_{11}(p), \hat \dd_{22}(1), \cdots,
\hat \dd_{\min\{m,n\}, \min\{m,n\}}(p) )\|_0 \le j } \left\{ \sum_{k=1}^p \sum_{i=1}^{\min \{ m,n \}} \left(\dd_{ii}(k) - \hat \dd_{ii}(k)\right)^2 \right\}.
\end{eqnarray*}
Here, $\| \cdot\|_0$ denotes the number of nonzero members of an argument. Thus, $$(\dd_{11}^{(j)}(1), \cdots, \dd_{11}^{(j)}(p), \dd_{22}^{(j)}(1), \cdots,
\dd^{(j)}_{\min\{m,n\}, \min\{m,n\}}(p) )$$ is an optimal solution of the last minimization problem in the above expressions.  This further implies that $\A_j \equiv \U *_L \S_j *_L \V^\top$ is an optimal solution of  $\min_{\B \in \BB_j} \|\A-\B\|_F^2$.  This completes the proof.
\qed

\section{Further Questions}

There are some further questions on this topic.

{\bf Question 1} What is the algebraic structure of the set of $p \times p$ doubly real-preserving unitary transformations? Example \ref{ex2.3} indicates that if $P$ is an orthogonal matrix and $F$ is the DFT matrix, $PF$ may not be real-preserving.  Thus, the set of $p \times p$ doubly real-preserving unitary transformations is not a group.  Proposition \ref{p3.1} partially identified the algebraic structure of this set.  However, it remains an interesting problem to fully characterize the algebraic structure of this set.   Similar questions can be raised for the set of  $p \times p$ real-preserving matrices, and the set of  $p \times p$ doubly real-preserving matrices.

{\bf Question 2} Can we extend eigenvalue analysis to tubal matrices in $\KK_p^{n \times n}$ for a doubly real-preserving unitary transformation?

{\bf Question 3} What part of matrix theory can be extended to tubal matrices?   What part of matrix theory can be extended to tubal matrices with substantial differences?   What part of matrix theory cannot be extended to tubal matrices?

\bigskip


\end{document}